\documentclass[11pt]{article}
\usepackage{amsmath}
\usepackage{amssymb}
\usepackage[usenames]{color}

 \usepackage{colortbl}
 \usepackage{mathrsfs}
\usepackage{mathrsfs}

\def\C{\mathbb{C}}

\newtheorem{theorem}{\hspace*{\parindent}Theorem}

\newtheorem{corollary_l}{\hspace*{\parindent}Corollary}

\newcounter{theremark}

\title{Two-point distortion theorems and the Schwarzian derivatives of meromorphic functions}
\author{V.N.\:Dubinin\footnote{E-mail:  \emph{dubinin@iam.dvo.ru}}
\\[10pt]\small{\textit{Far Eastern Federal University, Vladivostok, Russia}}
\\\small{\textit{Institute of Applied Mathematics, FEBRAS, Vladivostok, Russia}}}
\date{}
\begin{document}
\maketitle
\begin{abstract}
For a meromorphic function $f$ in the unit disk $U=\{z:\;|z|<1\}$ and arbitrary points $z_1,z_2$ in $U$ distinct from the poles of $f$, a sharp upper bound on the product $|f'(z_1)f'(z_2)|$ is established. Further, we prove a sharp distortion theorem involving the derivatives $f'(z_1)$, $f'(z_2)$ and the Schwarzian derivatives $S_f(z_1)$, $S_f(z_2)$ for $z_1,z_2\in U$. Both estimates hold true under some geometric restrictions on the image $f(U)$.
\end{abstract}

\bigskip


\bigskip


\section{Introduction}
Two-point distortion theorems for holomorphic and univalent functions have been investigated intensely over last two decades \cite[11-15]{8}.
A classic theorem of this type due to Goluzin asserts the following. For a holomorphic and univalent function $f$ in the unit disk $U=\{z:\;|z|<1\},$ and arbitrary points $z_1,\;z_2\in{U}$, put $w_k=f(z_k),\;k=1,2$. Then
\begin{equation}\label{1}
|(1-|z_1|^2)f'(z_1)(1-|z_2|^2)f'(z_2)|\tanh^2(d(z_1,z_2))\le |w_1-w_2|^2,
\end{equation}
where $d(z_1,z_2)$ denotes  the hyperbolic distance \cite[Chapter IV, Section 3, Theorem 1]{9} (see also \cite[p. 120]{14}). Equality in (1) is attained for two distinct points $z_1$, $z_2$ if $f$ maps the disk $U$ onto the complex plane $\overline{\C}_w$ cut along one or two rays lying on the straight line $L$ orthogonal to the line segment joining $f(z_1)$ and $f(z_2)$ and crossing it at the midpoint. If  $f$ is univalent then some estimates are also known involving both the derivatives $f'(z_1),\;f'(z_2)$ and the Schwarzian derivatives $S_f(z_1),\;S_f(z_2),\;z_1,z_2\in U$ \cite{1,5}, where
$$
S_f(z)=\left(\frac{f''(z)}{f'(z)}\right)'-\frac{1}{2}\left(\frac{f''(z)}{f'(z)}\right)^2.
$$
In the present paper we consider the following  problem:  how to relax the univalence property of $f$ and obtain two-point distortion theorems involving $f'$ and $S_f$ under some weaker assumptions. Chuaqui and Pommerenke  in  \cite{4} as well as Chuaqui et.al. in \cite{2} established the sharp bounds on the distortion $|f'(z_1)f'(z_2)|$ for the analytic and locally univalent functions in $U$ under certain restrictions imposed on the Schwarzian norm defined by
$$
||S_f||=\sup_{z\in U}(1-|z|^2)^2|S_f(z)|.
$$
An interesting work \cite{3} by Chuaqui and Osgood deals with other distortion theorems under the same restrictions (see also \cite{16}). Our approach here is different: we impose  geometric restrictions on the image $f(U)$.  Under such restrictions we establish the estimates of type (1) and inequalities involving the Schwarzian derivatives at two points distinct from the critical points and poles of a meromorphic function $f$.  In particular, Corollary 1 in Section 2 shows that inequality (1) holds for a meromorphic function $f$ defined in $U$ and any points $z_1,\;z_2\in{U}$ distinct from the poles of $f$ if every circular arc joining the points $f(z_1),\;f(z_2)$ that belongs to $f(U)$ is covered univalently by the function $f$\footnote[1]{In other words, every point of an arc that belongs entirely to $f(U)$ has single pre-image.}. More sophisticated geometric conditions on $f$ leading to a distortion theorem involving the Schwarzian derivatives are given in Section 3. Our approach goes back to \cite{6} (see also \cite{7}) and makes extensive use of the condenser approach elaborated in \cite{5}.

\section{Two-point distortion theorem}

In what follows, the Riemann surface $\mathscr{R}$  is viewed as a finite or countable number of plane domains glued together respecting the following conditions: the projection of each point of the surface $\mathscr{R}$ coincides with a point of a glued domain; a neighborhood of each point of $\mathscr{R}$ is either a univalent disk or a finite-sheeted disk with the unique ramification point located at its center (details of this model can be found in \cite[Part 3]{10}). When this may not cause confusion, we will not distinguish between the plane domains before gluing (which identifies some parts of the boundaries of these domains) and after doing so (when they become subdomains of $\mathscr{R}$). The  notions, statements, and notation from the book \cite{5} will be used with no special explanation. Most of them carry over to the Riemann surfaces of the type described above in a natural way. Let $\mathscr{R}$ be a Riemann surface over the $w$-plane and suppose that $W_1,W_2\in\mathscr{R}$ are distinct from the ramification points of $\mathscr{R}$.  Set $w_k ={\rm pr}W_k$,$k =1,2$ and assume that $w_1\neq w_2$, $w_k\neq\infty,\; k=1,2$.  Denote by $\Gamma(w_1,w_2)$ the collection of all circles in the $w$-plane passing through the points $w_1$ and $w_2$. Given  two domains $\mathscr{B}_1,\; \mathscr{B}_2$ on the surface $\mathscr{R}$, such  that   $W_k\in\mathscr{B}_k$, $k=1,2$, write $\mathscr{B}_k^*$, $k=1,2$,  for the set of all points in  $\mathscr{B}_k$ that can be joined with the point $W_k$ by a Jordan curve lying on the surface $\mathscr{R}$ not passing through the ramification  points and having the projection belonging to some circle from $\Gamma(w_1,w_2)$.  The next lemma was demonstrated in \cite{6}.

\medskip

{\bf Lemma} . {\it Suppose $({\rm pr}\mathscr{B}_1^*)\cap({\rm pr}\mathscr{B}_2^*)=\emptyset$. Then the following inequality holds{\rm :}                   \begin{equation}\label{2}r(\mathscr{B}_1,W_1)r(\mathscr{B}_2,W_2) \le |w_1-w_2|^2,
\end{equation}
where $r(\mathscr{B}_k,W_k)$ is the inner radius of the domain $\mathscr{B}_k$ with respect to the point $W_k,\;k=1,2$.}

\smallskip

For a meromorphic   function  $f$  defined in  the unit disk $U$, denote by $\mathscr{R}(f)$ the Riemann surface onto which $f$ maps $U$. A set $\Lambda$ on the surface $\mathscr{R}(f)$ is said to be one-sheeted if  different  points of $\Lambda$  have different projections onto the $w$-plane, i.e., if the condition $W',\; W''\in\Lambda, W'\neq W''$, implies  ${\rm pr}W'\not={\rm pr}W''$.
\begin{theorem}\hspace{-0.2cm}{\bf.}
Let $f:U\to\mathscr{R}(f)$ be a meromorphic function, and let $z_1,\;z_2$ be some points of $U$ satisfying $f'(z_k)\neq0,\;k=1,2,$ and such that the projections $w_k={\rm pr}f(z_k)$ meet the conditions $w_1\neq w_2$ and $w_k\neq\infty$, $k=1,2$.  Suppose $\gamma_1\cup\gamma_2$ is one-sheeted for any pair of non-intersecting Jordan curves  $\gamma_1,\gamma_2\in\mathscr{R}(f)$ not passing through the ramification points of $\mathscr{R}(f)$, lying over the same circle from $\Gamma(w_1,w_2)$ and such that $f(z_k)\in\gamma_k$, $k=1,2$. Then inequality {\rm (1)} holds and is sharp.
\end{theorem}
\textbf{Proof.}
We will first examine the case $z_1=-z_2=-\lambda$, where $0<\lambda< 1$. Let us introduce the following notation: $$G_1=\{z\in U:\;{\rm Re}z<0\},\;\;\;G_2=\{z\in U:\;{\rm Re}z>0\},$$ $$
\mathscr{B}_k=f(G_k),\;\;W_k=f(z_k),\;\;k=1,2.
$$
The invariance of the Green function under $f$ implies that
\begin{equation}\label{3}r(G_k,z_k)|f'(z_k)|=r(\mathscr{B}_k,W_k),\;\;k=1,2.
\end{equation}
We claim that
\begin{equation}\label{4}
({\rm pr}\mathscr{B}_1^*)\cap ({\rm pr}\mathscr{B}_2^*)=\emptyset,
\end{equation}
where $\mathscr{B}_1^*,\;\mathscr{B}_2^*$ have been defined before the Lemma. Indeed, if this condition is violated, then there exists a circle $\gamma\in\Gamma(w_1,w_2)$ and two Jordan curves $\gamma_1,\;\gamma_2$ on the surface $\mathscr{R}(f)$ lying over $\gamma$ and not passing through the ramification points such that
$$
({\rm pr}\mathscr{\gamma}_1)\cap ({\rm pr}\mathscr{\gamma}_2)\not=\emptyset\;\;{\rm and}\;\;W_k\in\gamma_k,\;k=1,2.
$$
Then $\gamma_1\cap\gamma_2\not=\emptyset$ by the hypotheses of the Theorem. However, this relation contradicts the fact that the domains $\mathscr{B}_1$ and $\mathscr{B}_2$ are non'overlapping which verifies that validity of (4). Hence, we are in the position to apply the Lemma which, combined with (3), yields the inequality
$$
\prod_{k=1}^2r(G_k,z_k)|f'(z_k)|\le |w_1-w_2|^2.
$$
A simple calculations shows that
$$
r(G_1,z_1)=r(G_2,z_2)=\frac{2\lambda(1-\lambda^2)}{1+\lambda^2}=(1-\lambda^2)\tanh(d(z_1,z_2)).
$$
Therefore, inequality (1) holds for $z_1=-\lambda,\;z_2=\lambda$.

The general case follows by an application of the above particular case to the composition $f\circ\varphi$, where $\varphi$ is the M$\ddot{\mbox{o}}$bius automorphism of $U$
such that $\varphi(-\lambda)=z_1,\;\varphi(\lambda)=z_2$ for some $\lambda\in(0,1)$ in view of  $d(-\lambda,\lambda)=d(z_1,z_2)$. This completes the proof of Theorem 1.

The formulation of the next corollary makes no explicit reference to the Riemann surface $\mathscr{R}(f)$.

\begin{corollary_l}\hspace{-0.2cm}{\bf.}
Let $f:\;U\to\overline{\C}_w$ be a meromorphic function, and let $z_1$, $z_2$ be points in $U$ distinct from the poles of $f$ and such that $f(z_1)\neq f(z_2).$ Suppose that every arc of the circle joining the points $f(z_1)$, $f(z_2)$ and belonging to the image $f(U)$ is univalently covered by $f$. Then inequality {\rm (1)} holds true with $w_k=f(z_k)$, $k=1,2$, and is sharp.
\end{corollary_l}

\section{An inequality involving the Schwarzian derivatives}

Given two distinct points $w_1,\;w_2\in\C_w $ and a real number $t$, $t\neq{0}$  define a pair of the closed Jordan curves by the equation
$$
\left|\left(\frac{2w-w_1-w_2}{w_2-w_1}\right)^2-1-i t\right|=|t|.
$$
These curves are symmetric to each other with respect to the point $(w_1+w_2)/2$. Denote by $\Delta(w_1,w_2)$ the family of these curves emerging when $t$ runs over the set $(-\infty,\infty)\setminus\{0\}$.
\begin{theorem}\hspace{-0.2cm}{\bf.}
Let $f:\;U\to\mathscr{R}(f)$ be a meromorphic function, and let $z_1,\;z_2$ be some points in $U$ such that $f'(z_k)\neq 0,\;k=1,2,$ and the projections $w_k={\rm pr}f(z_k)$ satisfy the conditions $w_1\neq w_2,\;w_k\neq\infty,\;k=1,2.$  Suppose that any Jordan curve on the surface $\mathscr{R}(f)$ passing through either of the points  $f(z_1)$, $f(z_2)$ but not through the ramification points of $\mathscr{R}(f)$ and lying over a curve from $\Delta(w_1,w_2)$ is one-sheeted.
Then the following sharp estimate holds{\rm :}
$$
{\rm Re}\left\{\sum_{k=1}^2\frac{S_f(z_k)(w_2-w_1)^2}{6(f'(z_k))^2}+\frac{2(w_2-w_1)^2}{f'(z_1)f'(z_2)(z_1-z_2)^2}+\frac{2|w_2-w_1|^2}{f'(z_1)\overline{f'(z_2)}(1-z_1\overline{z}_2)^2}\right\}\le$$
\begin{equation}\label{5}
\le 2+\frac{|w_2-w_1|^2}{|f'(z_1)|^2(1-|z_1|^2)^2}+\frac{|w_2-w_1|^2}{|f'(z_2)|^2(1-|z_2|^2)^2}.
\end{equation}
Equality in {\rm (5)} is attained, for example, for the functions of the form
$$
f(z)=\frac{z(1+\lambda^2)-i\lambda z^2-i\lambda}{\lambda z^2-i(1+\lambda^2)z+\lambda}
$$
and the points $z_1=-\lambda,\;z_2=\lambda,\;0<\lambda<1$.  Each extremal function $f$ maps $U$ conformally and univalently onto $w$-plane
slit along an arc of the circle $|w|=1$ maintaining the point correspondence $f(-\lambda)=-1$, $f(\lambda)=1$.
\end{theorem}
\textbf{Proof.}
We will first consider the case $f(z_2)=-f(z_1)=1.$ We may assume that the surface $\mathscr{R}(f)$ is bounded by an analytic curve. As $f'(z_1)\not=0$ and $f'(z_2)\not=0$, one can find non-overlapping open univalent disks $\mathscr{U}_1$ and $\mathscr{U}_2$ on $\mathscr{R}(f)$ centered at the points $f(z_1)$ and $f(z_2)$, respectively. Given $r,\rho>0$ introduce the notation
$$
\omega_1=-1-\rho,\;\;\omega_2=-1+\rho,\;\;\omega_3=1-\rho,\;\;\omega_4=1+\rho,
$$
$$
\Lambda_k=\{w:\;|w-\omega_k|\le r\},\;k=1,2,3,4.
$$
Choose $\rho>0$ sufficiently small to ensure the inclusions $\omega_1,\;\omega_2\in {\rm pr}\mathscr{U}_1,\;\omega_3,\;\omega_4\in {\rm pr}\mathscr{U}_2$. Next, consider the closed disks $\mathscr{E}_k(r)$ on the surface $\mathscr{R}(f)$ centered at the points $W_k$ such that
$$
{\rm pr}W_k=\omega_k,\;\;{\rm pr}\mathscr{E}_k(r)=\Lambda_k,\;k=1,2,3,4;
$$
$$
\mathscr{E}_1(r),\;\mathscr{E}_2(r)\subset \mathscr{U}_1,\;\;\mathscr{E}_3(r),\;\mathscr{E}_4(r)\subset \mathscr{U}_2.
$$
We assume that the disks $\Lambda_k,\;k=1,2,3,4$ are pairwise disjoint and do not intersect the imaginary axis.
It is always possible to secure this conditions as well as the above inclusions by taking $r>0$ small enough.

Define the condenser
$$
\mathscr{C}=(\mathscr{R}(f),\{\mathscr{E}_k(r)\}_{k=1}^4,\{\delta_k\}_{k=1}^4)
$$
on the surface $\mathscr{R}(f)$ with plates $\mathscr{E}_k(r)$, $k=1,2,3,4$,  and the potentials $\delta_1=-1,\;\delta_2=1,$ $\delta_3=1,\;\delta_4=-1$. Let $C$ be the  ''inverse image'' of the condenser $\mathscr{C}$ in the disk $U$.  More precisely,
$$
C=(U,\{E_k(r)\}_{k=1}^4,\{\delta_k\}_{k=1}^4),
$$
where $E_k(r)=f^{-1}(\mathscr{E}_k(r))$ is the closed ''almost disc'' of radius $r/|f'(\zeta_k)|$ centered at the point $\zeta_k:=f^{-1}(W_k)$ $k=1,2,3,4$. Using the conformal invariance of the condenser capacity and Theorem 2.1 from \cite{5} we obtain
$$
{\rm cap}\mathscr{C}={\rm cap}C=-\frac{8\pi}{\log r}-2\pi\left\{\sum_{k=1}^4\log[(1-|\zeta_k|^2)|f'(\zeta_k)|]+\right.
$$
\begin{equation}\label{6}
\left.+\sum_{k=1}^4\sum_{{\substack{
  l=1\\
 l\neq k}}}^4\delta_k\delta_l\log\left|\frac{1-\overline{\zeta}_k\zeta_l}{\zeta_k-\zeta_l}\right|\right\}\left(\frac{1}{\log r}\right)^2+o\left(\left(\frac{1}{\log r}\right)^2\right),\;\;r\to 0.
\end{equation}
Following \cite[Chapter 4.3]{5} we now define the separating transformation of the condenser $\mathscr{C}$ with respect to the family of functions $\zeta=p_k(w)\equiv(-1)^k i w^2,\;k=1,2,3,4.$ The function $\zeta=p_k(w)$ maps the sector $D_k:=\{w:\;\pi(k-1)/2<\arg w<\pi k/2\}$ conformally and univalently onto the right half-plane ${\rm Re}\zeta>0$. Introduce the notation:
$$
\mathscr{E}_{kj}(r)=\{W\in\mathscr{E}_j(r):\;{\rm pr}W\in\overline{D}_k\},\;\;k,j=1,2,3,4
$$
($\mathscr{E}_{kj}(r)$ may be an empty set for some values of $k$ and $j$). The set $\{W\in\mathscr{R}:{\rm pr}W\in{D_k}\}$ comprises a finite number of pairwise disjoint domains.  To each such domain we add its boundary points lying over $\partial{D}_k$ and denote the resulting collection of sets by $\mathscr{R}_k$. The function $p_k$ induces a mapping defined on each set from $\mathscr{R}_k$ onto some surface lying over ${\rm Re}\zeta\ge0$.   Denote by $\tilde{\mathscr{R}}_k$ the collection of these surfaces. We write $\tilde{\mathscr{E}}_{kj}(r)$ for the ''image''
of $\mathscr{E}_{kj}(r)$ under the mapping $p_k,\;k,j=1,2,3,4$. Let $\mathscr{H}_k$ be the union of the sets from  $\tilde{\mathscr{R}}_k$ with their reflections with respect to the imaginary axis, and let $\mathscr{E}^*_{kj}(r)$ be the union of the set $\tilde{\mathscr{E}}_{kj}(r)$ with its reflection with respect to the imaginary axis. The set  $\mathscr{H}_k$ is a union of some domains and $\mathscr{E}^*_{kj}(r)$ is a closed set in $\mathscr{H}_k$ (possibly empty). By the result of the separating transformation of the condenser $\mathscr{C}$ with respect to the family of functions $\{p_k(w)\}_{k=1}^4$ we mean the family $\{\mathscr{C}_k\}_{k=1}^4$ of the condensers
$$
\mathscr{C}_k=(\mathscr{H}_k,\{\mathscr{E}^*_{kj}(r)\}_{j},\{\delta_j\}_{j}),
$$
where the index $j$ runs over the values such that $\mathscr{E}^*_{kj}(r)\not=\emptyset$.  With the above notation we have
\begin{equation}\label{7}
{\rm cap}\mathscr{C}\ge \frac{1}{2}\sum_{k=1}^4{\rm cap}\mathscr{C}_k
\end{equation}
(the proof of inequality (7) is essentially the same as that of Theorem 4.8 in \cite{5}).

Next, we treat the case $k=1$ in detail. Non-empty sets among $\mathscr{E}^*_{1j},\;1\le j\le 4,$ are two closed ''almost discs'' centered at some points $Z_1,\;Z_2$ with ${\rm pr}Z_1=-i (1-\rho)^2,$ ${\rm pr}Z_2=-i (1+\rho)^2.$ The radius of the ''almost disc'' centered at $Z_1$ is equal to $2r(1-\rho)$, while the radius of the ''almost disc'' centered at $Z_2$ is equal to $2r(1+\rho)$. In accordance with \cite[Theorem 2.1]{5} we then have the following asymptotic formula:
$$
{\rm cap}\mathscr{C}_1=-\frac{4\pi}{\log r}-2\pi\left\{\log\frac{r(\mathscr{H}_1,Z_1)}{2(1-\rho)}+\log\frac{r(\mathscr{H}_1,Z_2)}{2(1+\rho)}-\right.
$$
\begin{equation}\label{8}
\left.-2g_{\mathscr{H}_1}(Z_1,Z_2)\right\}\left(\frac{1}{\log r}\right)^2+o\left(\left(\frac{1}{\log r}\right)^2\right),\;\;\;\;r\to 0,
\end{equation}
where $r(\mathscr{H}_1,Z_j)$ is the inner radius of the connected component of $\mathscr{H}_1$ containing the point $Z_j,$ and $g_{\mathscr{H}_1}(Z,Z_j)$ is the Green function of this component with pole at  $Z_j,\;j=1,2.$ Denote by $\mathscr{B}$ the connected component of $\mathscr{H}_1$ containing the points $Z_1,\;Z_2,$ and put
$$
\mathscr{B}_j=\{Z\in\mathscr{B}:\;(-1)^{j+1}(g_{\mathscr{B}}(Z,Z_1)-g_{\mathscr{B}}(Z,Z_2))>0\},\;\;\;j=1,2.
$$
The following equality holds:
\begin{equation}\label{9}
\log[r(\mathscr{H}_1,Z_1)r(\mathscr{H}_1,Z_2)]-2g_{\mathscr{H}_1}(Z_1,Z_2)=\log[r(\mathscr{B}_1,Z_1)r(\mathscr{B}_2,Z_2)].
\end{equation}
Indeed, the function $\mathscr{U}(Z):=g_{\mathscr{B}}(Z,Z_1)-g_{\mathscr{B}}(Z,Z_2)$ is harmonic in the set $\mathscr{B}\setminus\{Z_1,Z_2\}$ and  $\mathscr{U}(Z)\to +\infty$ as $Z\to Z_1$ while $\mathscr{U}(Z)\to -\infty$ as $Z\to Z_2.$ The Green function of the set $\mathscr{B}_1$ with pole at $Z_1$ coincides with $\mathscr{U}(Z)$ in $\mathscr{B}_1.$ Hence,
$$\log r(\mathscr{B}_1,Z_1)=\lim_{Z\to Z_1}(\mathscr{U}(Z)+\log|{\rm pr}Z-{\rm pr}Z_1|)=\log r(\mathscr{B},Z_1)-g_{\mathscr{B}}(Z_1,Z_2).$$
In a similar way,  the function $-\mathscr{U}(Z)$ coincides in $\mathscr{B}_2$ with the Green function of $\mathscr{B}_2$ with pole at $Z_2.$ Hence,
$$
\log r(\mathscr{B}_2,Z_2)=\log r(\mathscr{B},Z_2)-g_{\mathscr{B}}(Z_1,Z_2).
$$
Adding together the above relations and taking account of the equalities
$$
r(\mathscr{B},Z_j)=r(\mathscr{H}_1,Z_j),\;\;j=1,2,\;\;\;\;g_{\mathscr{B}}(Z_1,Z_2)=g_{\mathscr{H}_1}(Z_1,Z_2),
$$
we arrive at (9).

According to the hypothesis of the Theorem any Jordan curve on the surface $\mathscr{R}(f)$ passing through either of the points  $f(z_1)$, $f(z_2)$ but not through  the ramification points of $\mathscr{R}(f)$ that lies over a curve from $\Delta(w_1,w_2)$ is one-sheeted. It follows that any Jordan curve lying over a circle $\{\zeta:\;|\zeta+i-t|=|t|\},$ $0<|t|<\infty,$ that belongs to $\mathscr{H}_1$  and passes through the point $p_1(f(z_2))$\footnote[2]{Here, the point $f(z_2)$ is viewed as a point on the surface $\mathscr{R}(f)$ and the point $p_1(f(z_2))$ belongs to the set $\mathscr{H}_1$, ${\rm pr}(p_1(f(z_2)))=-i.$}, but not through the ramification points of $\mathscr{H}_1,$ is one-sheeted (a univalent curve in the terminology of \cite{6}). Repeating a part of the proof of the Theorem from pages 519-520 of \cite{6}, we deduce the inequality
\begin{equation}\label{10}
r(\mathscr{B}_1,Z_1)r(\mathscr{B}_2,Z_2)\le |{\rm pr}Z_1-{\rm pr}Z_2|^2=(4\rho)^2
\end{equation}
for sufficiently small $\rho>0$ (note that in \cite{6} the inner radii of the auxiliary domains are computed with respect to the points $W_1,\;W_2$, while here they are computed with respect to the points $Z_1,\;Z_2$). Summarizing the relations (8) -- (10) we obtain (for $k=1$)
\begin{equation}\label{11}{\rm cap}
\mathscr{C}_k\ge-\frac{4\pi}{\log r}-2\pi\left(\log\frac{4\rho^2}{1-\rho^2}\right)\left(\frac{1}{\log r}\right)^2+o\left(\left(\frac{1}{\log r}\right)^2\right),\;\;r\to 0.
\end{equation}
Similarly, inequality (11) can be established for $k=2,3$ and $4$. Combining this bound  with relations (6) and (7) we arrive at the resulting inequality
\begin{equation}\label{12}
\sum_{k=1}^4\log[(1-|\zeta_k|^2)|f'(\zeta_k)|]+\sum_{k=1}^4\sum_{{\substack{
  l=1\\
 l\neq k}}}^4\delta_k\delta_l\log\left|\frac{1-\overline{\zeta}_k\zeta_l}{\zeta_k-\zeta_l}\right|\le 2\log \frac{4\rho^2}{1-\rho^2}.
\end{equation}
Finally, we are interested in the limit case of (12) as $\rho\to0$. In some neighborhoods of the points $w=-1$ and $w=1$ there exists a holomorphic branch $h$ of the inverse function $f^{-1}$ such that $h(-1)=z_1,\;h(1)=z_2.$ In terms of the function $h$, inequality (12) takes the form
\begin{equation}\label{13}
\sum_{k=1}^4\log\frac{1-|h(\omega_k)|^2}{|h'(\omega_k)|}+\sum_{k=1}^4\sum_{{\substack{
  l=1\\
 l\neq k}}}^4\delta_k\delta_l\log\left|\frac{1-\overline{h(\omega_k)}h(\omega_l)}{h(\omega_k)-h(\omega_l)}\right|\le 2\log \frac{4\rho^2}{1-\rho^2}.
\end{equation}
We spell out the expansions of $h$ in a neighborhood of $-1$ and $1$ as follows:
$$
h(w)=z_1+a_1(w+1)+a_2(w+1)^2+a_3(w+1)^3+...,
$$
$$
h(w)=z_2+b_1(w-1)+b_2(w-1)^2+b_3(w-1)^3+...
$$
These formulas furnish the expansions of $h(\omega_k)$ and $h'(\omega_k)$ dependent on $\rho,\;k=1,2,3,4.$ By substituting these expansions into (13), performing some simple calculations and  letting $\rho\to0$ we obtain
$$
{\rm Re}\left\{-\frac{1}{6}S_{h}(-1)-\frac{1}{6}S_{h}(1)+\frac{2a_1b_1}{(z_1-z_2)^2}+\frac{2a_1\overline{b}_1}{(1-z_1\overline{z}_2)^2}\right\}\le
$$
$$
\le\frac{1}{2}+\frac{|a_1|^2}{(1-|z_1|^2)^2}+\frac{|b_1|^2}{(1-|z_2|^2)^2}.
$$
When rewritten in terms of $f$ the last inequality coincides with (5) for $f(z_2)=-f(z_1)=1$. In general situation, Theorem 2 is proved by an  application
of the particular case just established to the function
$$
\frac{2f-w_1-w_2}{w_2-w_1}.
$$
The equality case is straightforward to verify. This completes that proof of Theorem 2.

The following statement makes no mentioning of the Riemann surface $\mathscr{R}(f)$.
\begin{corollary_l}\hspace{-0.2cm}{\bf.}
Let $f:\;U\to\overline{\C}_w$ be a meromorphic function, and suppose that $z_1,\;z_2\in{U}$ are distinct from the poles and critical points of $f$ and have different images, $w_1=f(z_1)\neq w_2=f(z_2)$.  Assume further that every curve from the family $\Delta(w_1,w_2)$ that belongs to $f(U)$ is univalently covered by the function $f$.
Then inequality (5) holds true and is sharp.
\end{corollary_l}

\paragraph{Acknowledgements} This research was supported by the Russian Science Foundation, grant number 14-11-00022.

\end{document}